
\baselineskip=14pt
\parskip=10pt

\font\eightrm=cmr8 

\magnification=\magstephalf

\def\1{{\overline{1}}}
\def\2{{\overline{2}}}
\parindent=0pt
\overfullrule=0in

\def\frac#1#2{{#1 \over #2}}

\bf
\centerline
{
Explicit Expressions for the Variance and Higher Moments of the Size of 
}
\centerline
{
a Simultaneous Core Partition and its Limiting Distribution
}

\rm
\bigskip
\centerline
{\it By Shalosh B. EKHAD and Doron ZEILBERGER}
\bigskip
\qquad 
{\it 
Dedicated to William Y.C. ``Bill'' Chen, the tireless apostle of enumerative and algebraic combinatorics in China (and beyond)}
\bigskip

{\bf Important Update (Sept. 1, 2015)}: 
It turns out that the second challenge below (except for the page limit) has been met before our paper was written,
by Paul Johnson. See insightful comments by Marko Thiel and Nathan Williams:

{\tt http://www.math.rutgers.edu/\~{}zeilberg/mamarim/mamarimhtml/stcoreFeedback.html} \quad . 

Note, in particular, that all our of theorems are {\it now} {\bf rigorously proved} theorems. 
As mentioned by them, Theorems 2 and 3 have been anticipated in their paper 
``Strange Expectations''  

{\tt http://arxiv.org/abs/1502.07934} \quad .

A donation to the OEIS in honor of Paul Johnson, Marko Thiel, and Nathan Williams, has been made. 
The first challenge  remains wide open!

{\bf VERY IMPORTANT} 

As in all our joint papers, the main point is not the article, but the accompanying Maple package, {\tt stCore}, that
may be downloaded, free of charge, from the webpage of this article

{\tt http://www.math.rutgers.edu/\~{}zeilberg/mamarim/mamarimhtml/stcore.html} \quad ,

where the readers can also find sample input and output files, that they are welcome to extend using their
own computers.

{\bf Introduction}

Many (perhaps most) combinatorial {\it statistics} (e.g. the number of Heads in tossing a coin $n$ times, the number of
inversions of an $n$-permutation [and, more generally, the number of occurrences of  any pattern in an $n$-permutation]),
are {\it asymptotically normal}, which means that if you denote it by $X_n$, figure out the {\it average}, $a_n:={\bf E}[X_n]$,
and then figure out the {\it variance}, let's call it $m_2(n):={\bf E}[(X_n-a_n)^2]$, then the {\it centralized and standardized}
version $Z_n:=(X_n-a_n)/\sqrt{m_2(n)}$ tends, as $n \rightarrow \infty$, to the good old {\it normal distribution} 
(aka {\it Gaussian distribution})
whose probability density function is $\frac{1}{\sqrt{2\pi}}e^{-x^2/2}$. Our favorite way ([Z1][Z2]) of proving this is
{\it automatically}, by using {\it symbol-crunching} to compute (at least the leading terms) of
the general moment $m_r(n)$ and then prove, ({\it automatically}, of course) that 
$\lim_{n \rightarrow \infty} \frac{m_r(n)}{m_2(n)^{r/2}}$, equals $0$ for
$r$ odd, and $r!/((r/2)!\, 2^{r/2})$, for $r$ even, the famous moments of the normal 
distribution.

But there are numerous exceptions! The most notable is the statistics ``length of the longest increasing subsequence'' defined
over the set of $n$-permutations, where the intriguing Tracy-Widom distribution shows up. Another example is the
subject matter of the present article,  the random variable ``size'', defined over the set of $(s,t)$-core partitions
(see below for definition), where $s$ and $t$ are relatively prime positive integers.
One of us (DZ) pledges a $100$ dollars donation to the OEIS Foundation, in honor of the first prover(s), for meeting the following
challenge.

{\bf First Challenge}: Prove,  rigorously, that the scaled limiting distribution (see below) (as $(s,t)$ both go to infinity, and $s-t$ is a fixed constant) 
of the  combinatorial random variable ``size'' defined on the set of $(s,t)$-core partitions is given
by the continuous random variable
$$
\sum_{k=1}^{\infty} \frac{z_k^2 + \tilde{z}_k^2}{4 \pi^2 k^2} \quad ,
$$
where $z_k$ and $\tilde{z}_k$ are jointly independent sequences of independent standard normal random variables.

One of us (SBE) verified that the first
nine (standardized) moments of ``size of an $(s,t)$-core partition'', (as $(s,t) \rightarrow \infty$) 
converge to the corresponding moments of the above continuous distribution, and it is virtually certain that this is
true in general (see below for details).

This distribution is mentioned in [DGP], eq. (2.4), where it is called $U_{VS}(1)$.

But the main purpose of this article is to  show the power of {\it symbol-crunching} in deriving new mathematical knowledge.
We will state deep new (polynomial) expressions for  the first six moments of the random variable {\it size of an $(s,t)$-core partition},
and the first nine moments for the special case of $(s,s+1)$-core partitions. From the {\it ``religious-fanatical''}
viewpoint of the current ``mainstream'' mathematician, they are ``just'' conjectures, but nevertheless, they are {\bf absolutely certain}
(well, at least as absolutely certain as most proved theorems). We also briefly indicate how we derived
these expressions, and indicate, for those {\it obtuse} mathematicians who would like to see  mathematical proofs, how they may possibly be proved.

We are not offering any prizes for just any old proofs, but we will be delighted to donate another \$100 to the OEIS Foundation for
meeting the following challenge.

{\bf Second Challenge}: Come up with a ``hand-waving'' and ``soft'' (yet rigorous!) {\it a priori} reason 
(whose length is not to exceed two pages) why
the average, and any finite moment, must be a polynomial in $s$ and $t$, and also come up with a ``soft'' (but rigorous!)
upper bound for  the degrees. 

This would, {\it in one stroke}, prove, {\bf rigorously}, all the theorems in this article, because
it would rigorously justify the (empirically) obvious fact that they belong to the {\it polynomial ansatz},
and hence discoverable by {\it undetermined coefficients}, by gathering enough data, and  solving
a (usually large) system of linear equations, thereby turning the undetermined coefficients of the desired polynomial expressions 
into {\it determined} ones. (See below for details.)

Note that in many cases in combinatorics (e.g. [BZ], [E], [Z3]), such arguments are very easy, but in the present case, we don't see it.
We hope that one of our readers will!

{\bf (s,t)-Core Partitions and Drew Armstrong's Ex-Conjecture}

Recall that a {\it partition} is a non-increasing sequence of positive integers $\lambda=(\lambda_1, \dots, \lambda_k)$ with $k \geq 0$,
called its {\it number of parts}; $n:=\lambda_1 + \dots + \lambda_k$ is called its {\it size},
and we say that $\lambda$ is a {\it partition of $n$}.

Also recall that the {\it Ferrers diagram} 
(or equivalently, using empty squares rather than dots, {\it Young diagram}) 
of a partition $\lambda$ is obtained by placing, in a {\it left-justified} way,
$\lambda_i$ dots at the $i$-th row. For example, the Ferrers diagram of the partition $(5,4,2,1,1)$ is
$$
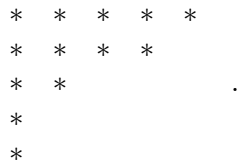
\matrix
{
* & * & *  &*  & * \cr
* & * & *  &*  &  \cr
* & * &   &  &  \cr
* &  &   &  &  \cr
* &  &   &  &  \cr
} \quad .
$$

Recall also that the {\it hook length} of a dot $(i,j)$ in the Ferrers diagram, $1\leq j \leq \lambda_i$, is the
number of dots to its right (in the same row) plus the number of dots below it (in the same column) plus one
(for itself), in other words $\lambda_i -i+\lambda'_j-j+1$, where $\lambda'$ is the {\it conjugate partition}, 
obtained by reversing the roles of rows and columns. (For example
if $\lambda=(5,4,2,1,1)$ as above, then $\lambda'=(5,3,2,2,1)$).

Here is a table of hook-lengths of the above partition, $(5,4,2,1,1)$:

$$
\matrix
{
9 & 6 & 4  & 3  & 1 \cr
7 & 4 & 1  &1  &  \cr
4 & 1 &   &  &  \cr
2 &  &   &  &  \cr
1 &  &   &  &  \cr
}  \quad .
$$
It follows that its set of hook-lengths is $\{1,2,3,4,6,7,9\}$. A partition is called an $s$-core if none of its hook-lengths is $s$.
For example, the above partition, $(5,4,2,1,1)$, is a $5$-core, and an $i$-core for all $i \geq 10$.

A partition is a {\it simultaneous} $(s,t)$-core partition if it avoids both $s$ and $t$. For example the
above   partition, $(5,4,2,1,1)$, is a $(5,11)$-core partition (and a $(5,12)$-core partition, and a $(100,103)$-core partition etc.).

For a lucid and engaging account, see [AHJ].

As mentioned in [AHJ], Jaclyn Anderson ([A]) very elegantly proved the following.

{\bf Theorem 0:} If $s$ and $t$ are relatively prime positive integers, then there are
{\bf exactly}
$$
\frac{(s+t-1)!}{s!t!} \quad,
$$
$(s,t)$-core partitions. 

For example, here are the $(3+5-1)!/(3!5!)=7$  $(3,5)$-core partitions:
$$
\{ empty , 1, 2, 11, 31, 211, 4211 \} \quad .
$$

Drew Armstrong ([AHJ], conjecture 2.6) conjectured, what is now the following theorem.

{\bf Theorem 1}: The {\it average size} of an $(s,t)$-core partition is given by
the nice polynomial
$$
\frac{(s-1)(t-1)(s+t+1)}{24} \quad .
$$

For example, the (respective) sizes of the above-mentioned $(3,5)$-core partitions are
$$
0,1,2,2,4,4,8 \quad,
$$
hence the average size is
$$
\frac{0+1+2+2+4+4+8}{7}=\frac{21}{7}=3 \quad,
$$
and this agrees with Armstrong's conjecture, since
$$
\frac{(3-1)(5-1)(3+5+1)}{24} \, = \, 3 \quad .
$$

Armstrong's conjecture was recently proved by Paul Johnson ([J]) using a very complicated (but ingenious!) 
argument (that does much more).
Shortly after, and almost {\it simultaneously} (no pun intended) it was re-proved by Victor Wang [Wan],
using another ingenious (and even more complicated) argument, that also does much more, in particular,
proving an intriguing conjecture of Tewodros Amdeberhan and Emily Sergel Leven ([AL]).
Prior to the full proofs by  Johnson and Wang, Richard Stanley and Fabrizio Zanello [SZ] came up with
a nice (but rather {\it ad hoc}) proof of the important special case of $(s,s+1)$-core partitions.

We should also mention the very interesting approach in [YZZ], that proved an important special case 
(also proved in [AL]) of
a conjecture of Amdeberhan and  Leven 
(but not directly related to the subject matter of the present article).

{\bf Higher Moments}

Recall that the $r$-th moment about the mean of a random variable, $X$,  is ${\bf E}[(X- {\bf E}[X])^r]$, where $\bf E$ is
the {\it expectation} operation.

{\bf Theorem 2}: If $s$ and $t$ are relatively prime positive integers, then
the {\it variance} (aka as the {\it second moment (about the mean)}) of the random variable ``size
of an $(s,t)$-core partition'', is given by the
nice polynomial expression:
$$
{\frac {1}{1440}}\,st \left( t-1 \right)  \left( s-1 \right)  \left( s+t+1 \right)  \left( s+t \right)  \quad .
$$
For example for $(3,5)$-core partitions it equals
$$
{\frac {1}{1440}}\,3 \cdot 5 \left( 5-1 \right)  \left( 3-1 \right)  \left( 3+5+1 \right)  \left( 3+5 \right) = 6 \quad,
$$
and indeed (since the average is $3$)
$$
\frac{(0-3)^2+(1-3)^2+ (2-3)^2 +(2-3)^2+ (4-3)^2+ (4-3)^2+ (8-3)^2}{7} =6 \quad .
$$

{\bf Theorem 3}: If $s$ and $t$ are relatively prime  positive integers, then
the {\it third moment (about the mean)} of the random variable ``size
of an $(s,t)$-core partition'' is given by the
nice polynomial expression:
$$
{\frac {1}{60480}}\,st \left( t-1 \right)  \left( s-1 \right)  \left( s+t+1 \right)  \left( s+t \right)  \left( 2\,{s}^{2}t+2\,s{t}^{2}-3\,{s}^{2}-3\,st-3\,
{t}^{2}-3 \right) \quad .
$$
For example for $(3,5)$-core partitions it equals
$$
{\frac {1}{60480}}\,3 \cdot 5 \cdot \left( 5-1 \right)  \cdot \left( 3-1 \right) \cdot  \left( 3+5+1 \right)  \cdot \left( 3+5 \right)  
\cdot \left( 2 \cdot {3}^{2} \cdot 5+2 \cdot 3 \cdot {5}^{2}-3 \cdot {3}^{2}-3 \cdot 3 \cdot 5-3 \cdot {5}^{2}-3 \right)= \frac{90}{7} \quad ,
$$
and indeed (since the average is $3$)
$$
\frac{(0-3)^3+(1-3)^3+ (2-3)^3 +(2-3)^3+ (4-3)^3+ (4-3)^3+ (8-3)^3}{7} = \frac{90}{7} \quad .
$$

{\bf Theorem 4}: If $s$ and $t$ are relatively prime  positive integers, then
the {\it fourth moment (about the mean)} of the random variable ``size
of an $(s,t)$-core partition'' is given by the
nice polynomial expression:
$$
{\frac {1}{4838400}}\,st \left( t-1 \right)  \left( s-1 \right)  \left( s+t+1 \right)  \left( s+t \right)  \cdot
$$
$$
( 19\,{s}^{4}{t}^{2}+38\,{s}^{3}{t}^{3}+19\,{s}^{2}{t}^{4}-51\,{s}^{4}t-102\,{s}^{3}{t}^{2}-102\,{s}^{2}{t}^{3}-51\,s{t}^{4}+36\,{s}^{4}
$$
$$
+72\,{s}^{3}t+108\,{s}
^{2}{t}^{2}+72\,s{t}^{3}+36\,{t}^{4}-33\,{s}^{2}t-33\,s{t}^{2}+36\,{s}^{2}+36\,st+36\,{t}^{2}+120 ) \quad .
$$
For example for $(3,5)$-core partitions it equals $\frac{726}{7}$,
and indeed (since the average is $3$)
$$
\frac{(0-3)^4+(1-3)^4+ (2-3)^4 +(2-3)^4+ (4-3)^4+ (4-3)^4+ (8-3)^4}{7} = \frac{726}{7} \quad .
$$

\vfill\eject

{\bf Theorem 5}: If $s$ and $t$ are relatively prime  positive integers, then
the {\it fifth moment (about the mean)} of the random variable ``size
of an $(s,t)$-core partition'' is given by the
nice polynomial expression (in computerish):

{\tt
1/95800320*s*t*(t-1)*(s-1)*(s+t+1)*(s+t)*(46*s**6*t**3+138*s**5*t**4+138*s**4*t**5+46
*s**3*t**6 
-211*s**6*t**2 - 633*s**5*t**3-844*s**4*t**4 - 633*s**3*t**5-211*s**2*t**6 + 333*s**6*t +
999*s**5*t**2 + 1665*s**4*t**3+1665*s**3*t**4 + 999*s**2*t**5+333*s*t**6 
- 180*s**6-540*s**5*t - 1283*s**4*t**2-1666*s**3*t**3 - 1283*s**2*t**4 - 540*s*t**5
-180*t**6 + 420*s**4*t + 840*s**3*t**2 + 840*s**2*t**3 + 420*s*t**4-180*s**4 - 360*s**3*t - 540*s**2*t**2
 - 360*s*t**3 - 180*t**4 + 327*s**2*t + 327*s*t**2 - 180*s**2 - 180*s*t - 180*t**2 - 3780) 
} \quad .

For example for $(3,5)$-core partitions it equals $\frac{2850}{7}$,
and indeed (since the average is $3$)
$$
\frac{(0-3)^5+(1-3)^5+ (2-3)^5 +(2-3)^5+ (4-3)^5+ (4-3)^5+ (8-3)^5}{7} = \frac{2850}{7} \quad .
$$

{\bf Theorem 6}: If $s$ and $t$ are relatively prime  positive integers, then
the {\it sixth moment (about the mean)} of the random variable ``size
of an $(s,t)$-core partition'' is given by the
nice polynomial expression (in computerish):

{\tt
1/4184557977600*s*t*(t-1)*(s-1)*(s+t+1)*(s+t)*(307561*s**8*t**4+1230244*s**7*t**5+
1845366*s**6*t**6+1230244*s**5*t**7+307561*s**4*t**8-2056306*s**8*t**3-8225224*s**7*t**4-
14394142*s**6*t**5-14394142*s**5*t**6-8225224*s**4*t**7-2056306*s**3*t**8+5372061*s**8*t
**2+21488244*s**7*t**3+42976488*s**6*t**4+53720610*s**5*t**5+42976488*s**4*t**6+21488244
*s**3*t**7+5372061*s**2*t**8-6453396*s**8*t-25813584*s**7*t**2-60704054*s**6*t**3-
91764618*s**5*t**4-91764618*s**4*t**5-60704054*s**3*t**6-25813584*s**2*t**7-6453396*s*t
**8+2985120*s**8+11940480*s**7*t+39743142*s**6*t**2+77437746*s**5*t**3+96285048*s**4*t**
4+77437746*s**3*t**5+39743142*s**2*t**6+11940480*s*t**7+2985120*t**8-11104272*s**6*t-
33312816*s**5*t**2-55521360*s**4*t**3-55521360*s**3*t**4-33312816*s**2*t**5-11104272*s*
t**6+2985120*s**6+8955360*s**5*t+23840061*s**4*t**2+32754522*s**3*t**3+23840061*s**2*t**
4+8955360*s*t**5+2985120*t**6-9109476*s**4*t-18218952*s**3*t**2-18218952*s**2*t**3-
9109476*s*t**4+2985120*s**4+5970240*s**3*t+8955360*s**2*t**2+5970240*s*t**3+2985120*t
**4+8664840*s**2*t+8664840*s*t**2-62687520*s**2-62687520*s*t-62687520*t**2+626875200) \quad .
}

For example for $(3,5)$-core partitions it equals $2346$,
and indeed (since the average is $3$)
$$
\frac{(0-3)^6+(1-3)^6+ (2-3)^6 +(2-3)^6+ (4-3)^6+ (4-3)^6+ (8-3)^6}{7} = \frac{16422}{7}=2346 \quad .
$$

The last three theorems regard the special case of $(s,s+1)$-core partitions.

{\bf Theorem 7}: If $s$ is a positive  integer, then
the {\it seventh moment (about the mean)} of the random variable ``size
of an $(s,s+1)$-core partition'' is given by the
nice polynomial expression (in computerish):

{\tt
1/149448499200*s**2*(s-1)*(s-2)*(2*s+1)*
(124496*s**14 - 527660*s**13 - 127268*s**12+ 2133077*s**11
+ 1565655*s**10-3928575*s**9-7848989*s**8-3573289*s**7
+ 7257797*s**6 +16741975*s**5+16528197*s**4+3583272*s**3
- 67819248*s**2-18541440*s+138620160) \hfill\break
*(s+1)**2
} \quad .

{\bf Theorem 8}: If $s$ is a positive  integer, then
the {\it $8$-th moment (about the mean)} of the random variable ``size
of an $(s,s+1)$-core partition'' is given by the
nice polynomial expression (in computerish):

{\tt
1/914624815104000*s**2 * (s-1)*(2*s+1)* 
(308851624*s**18 - 2759073420*s**17 + 7345195650*s**16 + 1614779679*s**15 - 27716691813*s**14
- 3203324556*s**13 + 61922226136*s**12 +
52270343442*s**11 - 49025878614*s**10 - 146716496688*s**9
- 153171599682*s**8 - 30342055161*s**7
+ 158893451131*s**6 - 165853921776*s**5 + 1073038790016*s**4
+ 9260929255680*s**3 - 11293714925568*s**2 - 19188060088320*s + 21924617379840)* (s+1)**2
} \quad .

{\bf Theorem 9}: If $s$ is a positive  integer, then
the {\it ninth moment (about the mean)} of the random variable ``size
of an $(s,s+1)$-core partition'' is given by the
nice polynomial expression (in computerish):

{\tt
1/182467650613248000*s**2*(s-1)*(s-2)*(2*s+1)*
(28092743584*s**20 - 284614603048*s**19 + 908242721124*s**18 - 87722680542*s**17
- 4040707469643*s**16 + 1347179583168*s**15 + 11350317109273*s**14
+ 4824122583716*s**13 - 15816684214230*s**12 - 31535118689736*s**11
- 29475404073738*s**10 + 2671156715274*s**9 + 63014451511513*s**8
 + 79700408583680*s**7 + 45859575725901*s**6 - 377516262865248*s**5
 + 6309067352294376*s**4 + 10737857697068736*s**3 - 38301852570773760*s**2
- 26103018295756800*s + 48704747653094400 )* (s+1)**2
} \quad .

{\bf A Crash course in Combinatorial Statistics}

Recall ([Z1][Z2]) that given a sequence of combinatorial random variables 
(e.g. the number of Heads upon tossing a fair coin $n$ times),
whose combinatorial generating function is either known explicitly 
($C_n(t)=(1+t)^n$ in this trivial  case), or only as expressions in $t$, but for many $n$,
one first finds the discrete {\it probability generating function} $P_n(t):=C_n(t)/C_n(1)$
(e.g. $P_n(t)=(1+t)^n/2^n$ for coin-tossing),
(under the uniform distribution). Of course $P_n(1)=1$ as it should.
To find the {\it first moment}, $av_n$ (alias average, alias mean, alias expectation) one 
either derives explicitly $P'_n(1)$ ($n/2$ in this trivial case), or ``guesses'' in more complicated cases.
Next one {\it centralizes} getting the centralized probability generating function
$\tilde{P}_n(t)=P_n(t)/t^{av_n}$, and gets the variance (alias second  moment about the mean),
$m_2(n)$, by computing (or ``guessing'') $(t \frac{d}{dt})^2 \tilde{P}_n(t)|_{t=1}$, and more generally,
the higher moments, $m_k(n)$, are given by $(t \frac{d}{dt})^k \tilde{P}_n(t)|_{t=1}$. From these
one derives the {\it standardized} (scaled) moments
$\alpha_k(n):=\frac{m_k(n)}{m_2(n)^{k/2}}$. It is often the case,
as $n \rightarrow \infty$, that the $\alpha_k(n)$ tend to fixed numbers, in which case we have
a {\it limiting distribution}. 

The above is easily extended to the case where the sequence of combinatorial random variables depends on
several discrete parameters, like the present case where they depend on the two parameters $s$ and $t$.

{\bf The Wonder of Internet Searches: How we found the (so far conjectured, but absolutely certain) Scaled Limiting Distribution}

As mentioned above, more often than not, this limiting
distribution is the good-old Gaussian, but this is {\it definitely} not the case this time.
Using our Theorems 2-4, we computed that the limiting {\it skewness}  (standardized third moment) happens to be
$$
\lim_{s,t \rightarrow \infty} \alpha_3(s,t) = \frac{4}{7} \cdot \sqrt{10}  \approx 1.807  \quad , 
$$
while the limiting {\it kurtosis} (standardized fourth moment) turns out to be
$$
\lim_{s,t \rightarrow \infty} \alpha_4(s,t) = \frac{57}{7} \approx 8.1429  \quad .
$$

Being citizens of our time, we {\it googled}

{\tt 1.807 8.1429 statistics}

and, {\it lo and behold}, got  the link

{\tt www.aueb.gr/conferences/Crete2015/Papers/Dalla.pdf} \quad ,

(reference [DGP]), that mentioned that these are the skewness and kurtosis of the 
continuous probability distribution that they call $U_{VS}(1)$ (already mentioned in the introduction), 
but we will  call $Z$
$$
Z \, := \,
\sum_{k=1}^{\infty} \frac{z_k^2 + \tilde{z}_k^2}{4 \pi^2 k^2} \quad ,
$$
where $z_k$ and $\tilde{z}_k$ are jointly independent sequences of independent normal random variables.

Liudas Giraitis, one of the authors of [DGP], kindly offered the following interesting information via email.

{\eightrm
``VS''  stands for ``rescaled  variance statistic".
Characterization of  the $U_{VS}(1)$   distribution as  the  sum of the  weighted sum of iid  normals
was well- known in statistical literature.
It  was established  by Geoffrey S. Watson ([Wat]) in the context of goodness-of-fit tests on a circle,
see  formula  (15) in his  paper.  
Later  in [GKLT], the authors 
show  that the   VS statistic  has  the  same   limit  as  the  limit  derived  by  
Watson  in the context of goodness-of-fit tests on a circle.
}

{\bf A Human Intermezzo: A Crash Course in Moment Generating Functions}

This inspired us to compute higher moments of the continuous probability distribution, $Z$,  as follows.

Recall that the {\it moment generating function} of a probability distribution $X$ is
the (exponential) generating function
$$
M_X(t) := \sum_{k=0}^{\infty} \frac{m_k}{k!} t^k \quad .
$$
It is well-known and trivial to see that, for any fixed constant, $c$,
$$
M_{cX}(t)=M_X(ct) \quad,
$$
and it is also well-known and easy (but not utterly trivial) to see that if $X$ and $Y$ are {\it independent} random variables, then
$$
M_{X+Y}(t)=M_X(t) M_Y(t) \quad .
$$

Hence if $\{X_i\}_{i=1}^{\infty}$ is a sequence of pairwise independent random variables and $\{c_i\}_{i=1}^{\infty}$ is a sequence
of positive numbers such that their sum converges, then
$$
M_{\sum_{i=1}^{\infty} c_i X_i} (t) =\prod_{i=1}^{\infty} M_{X_i}(c_i t) \quad .
$$

If, in addition, the $X_i$'s are {\it identically distributed}, denoting $M_{X_i}(t)$ by $M(t)$, we have
$$
M_{\sum_{i=1}^{\infty} c_i X_i} (t) =\prod_{i=1}^{\infty} M (c_i t) \quad .
$$

The famous moments of the standard normal distribution, $z$, are, as mentioned above, $0$ for $r$ odd and $r!/(2^{r/2} (r/2)!)$ for
$r$ even. Hence the $r$-th moment of $z^2$ is $(2r)!/(2^r r!)$. Hence the (exponential) moment generating function is
$$
M_{z^2}(t) \, = \,
\sum_{r=0}^{\infty} \frac{(2r)!}{2^r r!} \cdot \frac{t^r}{r!} =
\sum_{r=0}^{\infty} \frac{(2r)!}{2^r r!^2} t^r =
(1-2t)^{-1/2} \quad.
$$

Hence the (exponential) moment generating function of $z^2+\tilde{z}^2$, where $z$ and $\tilde{z}$ are independent
standard normal distributions, is
$$
M_{z^2+ \tilde{z}^2}(t)= \left( (1-2t)^{-1/2} \right)^2 \, = \, (1-2t)^{-1} \quad.
$$
Hence the (exponential) moment generating function of the continuous distribution $Z$ is
$$
M_Z(t) =\prod_{k=1}^{\infty} (1-\frac{t}{2 \pi^2 k^2})^{-1} \quad .
$$
But thanks to good old Leonhard's {\it iconic} $\sin(\pi x)= \pi x \prod_{k=1}^{\infty} (1- \frac{x^2}{k^2})$,
this equals
$$
\frac{\sqrt{t/2}}{\sin \sqrt{t/2}} \quad ,
$$
and this can be used to find as many (straight) moments as desired (take the coefficient of $t^k$ in the Maclaurin expansion
and multiply by $k!$). From these straight moments one easily computes the {\it moments about the mean}, using, thanks
to the binomial theorem,
${\bf E}[(X-m_1)^k]=\sum_{i=0}^{k} {{k} \choose {i}} (-m_1)^{k-i} m_i$, and from them the standardized moments.

[This is implemented in procedure {\tt VSmoms(N)} in the Maple package {\tt stCore}].

Using this, one of us (SBE) found:
$$
\alpha_3= \frac{4}{7} \cdot \sqrt{10} \quad ,  \quad
\alpha_4= \frac{57}{7} \quad , \quad
\alpha_5= \frac{820}{77} \cdot \sqrt{10} \quad , \quad
$$
$$
\alpha_6=  \frac{1537805}{7007} \quad, \quad
\alpha_7=  \frac{466860}{1001} \cdot \sqrt{10} \quad, \quad
\alpha_8=  \frac{193032265}{17017} \quad, \quad
$$
$$
\alpha_9 =\frac{70231858960}{2263261} \cdot \sqrt{10} \quad ,
$$

and these coincide {\it exactly} (up to the ninth moment) with the limiting (scaled) moments of our 
combinatorial random variable ``size of an $(s,t)$-core partition'', as $(s,t)$ go to infinity,
(and $s-t$ is bounded).

{\bf How did we derive the above Explicit Expressions for the First Nine Moments?}

Let us now briefly describe how we were able to discover the above very deep theorems regarding
the first nine moments.

The first step was to use Anderson's bijection, as modified in [AHJ], and to work with
$(s,t)$-{\bf Dyck Paths}, that are two-dimensional walks, using horizontal and vertical positive unit steps,
from the origin to the point $(s,t)$ staying {\it above} the line $s\,y\, - \,t\,x \, = \, 0$. 

[{\eightrm Equivalently, the number of ways of tossing a coin $t+s$ times such that you win $s$ dollars if it lands
Heads and lose $t$ dollars if it is Tails, breaking even at the end, but never being in debt}].

Each such path has associated with it a certain set of positive integers, 
nicely described in [AHJ], and also important is
how many such positive integers there are. Introducing the formal variables
$q$ and $w$, we associate the {\it weight} $q^s w^k$ to each such path,
where $s$ is the sum of the labels and $k$ is the number of such positive labels
(see [AHJ] for exactly how to determine them).

The same definitions apply to {\it partial walks} from $(0,0)$ to $(i,j)$ where
we look at the vertical steps $(i,j-1) \rightarrow (i,j)$, and we naturally 
associate the partial weight due to each such vertical step 
(assuming $(s,t)$ are fixed) as follows (for the sake of typographical clarity we use $x**y$ for $x^y$).
$$
Wt(i,j)(q,w):=
q**\left( \sum_{\{i' \geq i \, ; \,  s\,j-t\,i'-b>0\}} s\,j-t\,i' -b \right) \cdot w**\left( \sum_{\{i'\geq i \, ; \,  s\,j\, - \, t\,i'-b>0 \}} 1 \right) \quad ,
$$
then the {\it weight-enumerator} of the set of partial walks from $(0,0)$ to $(i,j)$ 
may be computed by the recurrence (we suppress the implied dependence on $(q,w)$)
$$
F_{s,t}(i,j) \, = \, F_{s,t}(i-1,j)+ Wt(i,j-1) \, \cdot \, F_{s,t}(i,j-1) \quad .
$$

We impose the obvious {\it initial condition}  $F_{s,t}(0,0)=1$, and {\it boundary conditions}
$F_{s,t}(i,j)=0$  if $i<0$ or $s\,j\, - \,t\,i \, < \, 0$.

Using this {\it linear recurrence scheme}, one can compute these weight-enumerators very fast.

But $F_{s,t}(i,j)(q,w)$ is but a {\it stepping stone}. At the end of the day
we take $F_{s,t}(s,t)(q, w)$, expand it in terms of powers of w and perform the
{\it umbral substitution}
$$
w^k \rightarrow q^{-k(k-1)/2} \quad .
$$
(This corresponds to going from the vector of hook-lengths of the first column of
the $(s,t)$-core partition to its shape, see [AHJ]).

All this is implemented in procedure {\tt Fab(a,b,q)} in our Maple package {\tt stCore}.

Using this efficient, {\it Dynamical Programming}, approach we can
{\it crank out} these weight-enumerators, $F_{s,t}(q)$, for many choices of coprime pairs $(s,t)$, 
and by repeatedly applying the operator $q \frac{d}{dq}$, and then plugging-in $q=1$
collect numerical data for the respective moments, and from these, numerical data (for many distinct $(s,t)$)
about the moments about the mean. 
Then, inspired by Drew Armstrong's ex-conjecture, that expressed the average as a {\bf polynomial},
we assume  (without any {\it a priori} theoretical justification!, see the second  challenge)
that these are polynomials (i.e. we use the {\it polynomial ansatz}),
and use {\it undetermined coefficients}, and fit the data with polynomials.

{\bf A possible Approach for (most probably ugly) mathematical proofs}

The difficulty with the $F_{s,t}(q)$, in their initial definition as weight-enumerators, according to ``size
of $(s,t)$-core partitions'', is that they are only defined for $s$ and $t$ relatively prime.
While the notion of $(s,t)$-Dyck path makes sense when $s$ and $t$ are not relatively prime, 
the nice formula that enumerates them, $(s+t-1)!/(s!t!)$, is no longer valid.

A natural approach would be not to be hung up on ``nice'' formulas, and set-up a general recursion schemes
for the weight-enumerators of partial walks, $F_{s,t}(i,j)(q,w)$, but no longer restricted to relatively
prime $s,t$. Then get some recurrence scheme that would enable an (ugly, ad-hoc)  proof in the
so-called {\it holonomic ansatz}. But since this approach is only likely to produce ugly proofs, and
we are {\it absolutely certain} that our Theorems 2-9 are correct, we rather not even try,
and do not encourage anyone to follow this approach.
(But those who insist should start with the easier special case when $t=s+1$.)

{\bf But we will be {\it thrilled} if any of our readers would meet our two challenges stated in the introduction!}

{\bf Acknowledgment}

We wish to thank Tewodros Amdeberhan for introducing us to this interesting subject.
Also many thanks to Liudas Giraitis for  explaining the motivation for the interesting paper [DGP],
and for pointing us to references [Wat] and [GKLT].

This project was started while we were visiting the Center of Combinatorics at Nankai
University, Tianjin, China, and benefited from the great hospitality of, and numerous stimulating
conversations with, the founding director, William  Y.C. ``Bill'' Chen, (to whom this paper is dedicated), 
the deputy director Arthur Yang, and Professor Qing-Hu Hou.

{\bf References}

[AL] Tewodros Amdeberhan and Emily Leven, {\it Multi-cores, posets, and lattice paths}, \hfill\break
{\tt http://arxiv.org/abs/1406.2250} \quad .

[A] Jaclyn Anderson,  {\it Partitions which are simultaneously $t_1$ - and $t_2$ -core}. Disc. Math. {\bf 248} (2002),  237-243.

[AHJ] Drew Armstrong, Christopher R. H. Hanusa, and Brant C. Jones,
{\it Results and conjectures on simultaneous core partitions}, \hfill\break
{\tt http://arxiv.org/abs/1308.0572} \quad .

[BZ] Andrew Baxter and Doron Zeilberger,
{\it The Number of Inversions and the Major Index of Permutations are Asymptotically Joint-Independently Normal}, 
the Personal Journal of Shalosh B. Ekhad and Doron Zeilberger, posted Feb. 4, 2011,
\hfill\break
{\tt http://www.math.rutgers.edu/\~{}zeilberg/mamarim/mamarimhtml/invmaj.html} \quad .

[DGP] Violetta Dalla , Liudas Giraitis, and Peter C.B. Phillips, {\it Testing Mean Stability of a Heteroskedastic Time Series}, preprint,
\hfill\break
{\tt www.aueb.gr/conferences/Crete2015/Papers/Dalla.pdf} \quad .

[E] Shalosh B. Ekhad, {\it The joy of brute force: the covariance of the number of inversions and the major index},
the Personal Journal of Shalosh B. Ekhad and Doron Zeilberger, posted ca. June 1995,
\hfill\break
{\tt http://www.math.rutgers.edu/\~{}zeilberg/mamarim/mamarimhtml/brute.html} \quad .

[GKLT] L. Giraitis, P. Kokoszka, R. Leipus, and G. Teyssiere,
{\it Rescaled variance and related tests for long memory in volatility and levels}, Journal of Econometrics {\bf 112}  (2003), 265-294.

[J] Paul Johnson, {\it Lattice points and simultaneous core partitions}, \hfill\break
{\tt http://arxiv.org/abs/1502.07934} \quad .

[SZ] Richard P. Stanley and Fabrizio Zanello, {\it The Catalan case of Armstrong's conjecture on simultaneous core partitions},  \hfill\break
{\tt http://arxiv.org/abs/1312.4352} \quad .

[Wan] Victor Y. Wang, {\it Simultaneous core partitions: parameterizations and sums}, \hfill\break
{\tt http://arxiv.org/abs/1507.04290} \quad .

[Wat] G.S. Watson, {\it Goodness-Of-Fit Tests on a Circle}, Biometrika, {\bf 48} (1961), 109-114.

[YZZ] Jane Y.X. Yang, Michael X.X. Zhong, and Robin D.P. Zhou, {\it On the Enumeration of (s,s+1,s+2)-Core Partitions}, \hfill\break
{\tt http://arxiv.org/abs/1406.2583v2} \quad .

[Z1] Doron Zeilberger, {\it The Automatic Central Limit Theorems Generator (and Much More!)},
in: ``Advances in Combinatorial Mathematics: Proceedings of the Waterloo Workshop in 
Computer Algebra 2008 in honor of Georgy P. Egorychev'', chapter 8, pp. 165-174, (I.Kotsireas, E.Zima, eds. Springer Verlag, 2009),
\hfill\break
{\tt http://www.math.rutgers.edu/\~{}zeilberg/mamarim/mamarimhtml/georgy.html} \quad .

[Z2]  Doron Zeilberger, {\it HISTABRUT: A Maple Package for Symbol-Crunching in Probability theory},
the Personal Journal of Shalosh B. Ekhad and Doron Zeilberger, posted Aug. 25, 2010,
\hfill\break
{\tt http://www.math.rutgers.edu/\~{}zeilberg/mamarim/mamarimhtml/histabrut.html} \quad .

[Z3]  Doron Zeilberger, {\it Symbolic Moment Calculus I.: Foundations and Permutation Pattern Statistics},
Annals of Combinatorics {\bf 8} (2004), 369-378,
\hfill\break
{\tt http://www.math.rutgers.edu/\~{}zeilberg/mamarim/mamarimhtml/smcI.html} \quad .

\bigskip
\bigskip
\hrule
\bigskip
Doron Zeilberger, Department of Mathematics, Rutgers University (New Brunswick), Hill Center-Busch Campus, 110 Frelinghuysen
Rd., Piscataway, NJ 08854-8019, USA. \hfill \break
zeilberg at math dot rutgers dot edu \quad ;  \quad {\tt http://www.math.rutgers.edu/\~{}zeilberg/} \quad .
\bigskip
\hrule
\bigskip
Shalosh B. Ekhad, c/o D. Zeilberger, Department of Mathematics, Rutgers University (New Brunswick), Hill Center-Busch Campus, 110 Frelinghuysen
Rd., Piscataway, NJ 08854-8019, USA.
\bigskip
\hrule

\bigskip
Exclusively published in The Personal Journal of Shalosh B. Ekhad and Doron Zeilberger  \hfill \break
({ \tt http://www.math.rutgers.edu/\~{}zeilberg/pj.html})
and {\tt arxiv.org} \quad . 
\bigskip
\hrule
\bigskip
{\bf First Version: Aug. 30, 2015} \quad; \quad {\bf This Version: Sept. 1, 2015} \quad; \quad

\end